\documentclass[12pt]{amsart}

\usepackage{amsfonts,amssymb,amscd}
\newcommand{\Q}{{\mathbb Q}}

\newcommand{\ti}{\tilde}

\newcommand{\odd}{{\mathrm{Odd}\,}}
\newcommand{\even}{{\mathrm{Even}\,}}
\newcommand{\ad}{{\mathrm{ad}\,}}

\newcommand{\lra}{\leftrightarrow}

\newtheorem{theorem}{Theorem}[section]

\newtheorem{proposition}[theorem]{Proposition}
\newtheorem{lemma}[theorem]{Lemma}
\newtheorem{claim}[theorem]{Claim}

\newtheorem{definition}[theorem]{Definition}
\newtheorem{example}[theorem]{Example}

\title[The Baker-Campbell-Hausdorff formula]
{The Baker-Campbell-Hausdorff formula \\
 in the free metabelian Lie algebra}

\author[V.~Kurlin]{V.~Kurlin}
\address{
 Division of Pure Mathematics,
 Department of Mathematical Sciences,
 University of Liverpool,
 Liverpool L69 7ZL, United Kingdom}
\email{ kurlin@liv.ac.uk }


\date{ December 10, 2006, the last version is at www.geocities.com/vak26}

\begin{document}

\begin{abstract}
The classical Baker-Campbell-Hausdorff formula gives a recursive way
 to compute the Hausdorff series $H=\ln(e^Xe^Y)$ for
 non-commuting $X,Y$.
Formally $H$ lives in the graded completion of
 the free Lie algebra $L$ generated by $X,Y$.
We present a closed explicit formula for $H=\ln(e^Xe^Y)$
 in a linear basis of the graded completion of 
 the free metabelian Lie algebra $\bar L=L/[[L,L],[L,L]]$.
\end{abstract}

\maketitle


\section{Introduction}


\subsection{Brief summary}
\noindent
\smallskip

The Baker-Campbell-Hausdorff (BCH) formula arises naturally
 in the context of Lie groups and Lie algebras.
Originally the series $H=\ln(e^Xe^Y)$
 was used to define a multiplication law in a Lie group
 associated to a given Lie algebra.
If the variables $X,Y$ commute then $\ln(e^Xe^Y)=X+Y$.
\smallskip

Let $L$ be the free Lie algebra generated by $X,Y$.
Then $H=\ln(e^Xe^Y)$ belongs to the graded completion of $L$, 
 i.e. $H$ contains commutators of all degrees (lengths).
A linear basis of $L$ contains exponentially many
 elements of a fixed degree. 
That is why the classical BCH formula
 is awkward for solving exponential
 equations in Lie algebras.
\smallskip

Author's interest in the BCH formula
 came from knot theory and number theory.
The Kontsevich integral is a powerful knot invariant
 and can be computed combinatorially from a knot projection 
 via a Drinfeld associator \cite{Dr1}.
A Drinfeld associator is a non-commutative 2-variable series
 living in the graded completion of a Lie algebra and 
 satisfying the pentagon and hexagon equations 
 involving 5 and 6 exponential factors, respectively.
An outstanding problem in quantum algebra is to compute explicitly
 a Drinfeld associator with rational coefficients.
\smallskip

The set of Drinfeld associators modulo commutators of commutators 
 contains a specific transcedental solution 
 expressed via the classical zeta values
 $\zeta(n)=\sum\limits_{k=1}^{\infty}\dfrac{1}{k^n}$.
A closed metabelian version of the BCH formula was used 
 to solve completely the pentagon and hexagon equations 
 modulo commutators of commutators \cite[Theorem~1.5c]{Kur}.
It turned out that these equations do not contain 
 polynomial relations between odd zeta values.
Here we present an elementary proof of the explicit BCH formula 
 for $H=\ln(e^Xe^Y)$ in the free metabelian Lie algebra 
 $\bar L=L/\bigl[[L,L],[L,L]\bigr]$.
\smallskip


\subsection{Definitions and results}
\noindent
\smallskip

Classical notions of the Lie theory can be found in \cite{Reu}.
Lie algebras are considered over a field of characteristic 0.
The free Lie algebra $L$ generated by $X,Y$
 is graded by the \emph{degree} : $\deg X=\deg Y=1$ and 
 $\deg[A,B]=\deg A+\deg B$ for all $A,B\in L$.
The \emph{graded} completion $\hat L$ of the Lie algebra $L$
 is the algebra of infinite series of elements of $L$.
\smallskip

The \emph{Hausdorff} series is $H=\ln(e^Xe^Y)$, where
 the logarithm and exponential are considered as formal power series, 
 i.e. $e^X=\sum\limits_{n=0}^{\infty}\dfrac{X^n}{n!}$
 and $\ln(1+Y)=\sum\limits_{n=1}^{\infty}\dfrac{(-1)^{n-1}}{n}Y^n$.
The classical BCH formula states that $H\in\hat L$ \cite{Bak, Cam, Hau}, 
 i.e. $H$ can be expressed as an infinite sum of commutators,
 see Theorem~2.4.
For $A_1,\ldots,A_m\in L$, the \emph{long commutator}
 $[A_1A_2\ldots A_{m-1}A_m]$ is `the right bracketing'
 $[A_1,[A_2,[\ldots[A_{m-1},A_m]\ldots]]]]$,
 e.g. $[X^2Y]=[X,[X,Y]]$ and $[YXY]=[Y,[X,Y]]$.
Then
$$H=X+Y+\dfrac{[XY]}{2}+\frac{[X^2Y]-[YXY]}{12}-\frac{[XYXY]}{24}+\cdots$$

E.~Dynkin found a closed formula for $H$ \cite{Dyn}, but not in 
 a linear basis of the graded completion $\hat L$, see Theorem~2.5.
Also the series $H$ can be expressed via associative monomials $W$
 in the variables $X,Y$ as follows: $H=X+Y+\sum c_W W$.
The generating function for the coefficients $c_W$
 was computed explicitly by K.~Goldberg \cite{Gol}, see Theorem~2.6.
\smallskip

If a Lie algebra $L$ satisfies $\bigl[[L,L],[L,L]\bigr]=0$ 
 then $L$ is said to be \emph{metabelian}.
For a free Lie algebra $L$, the quotient 
 $\bar L=L/\bigl[[L,L],[L,L]\bigr]$ is sometimes 
 called the \emph{free metabelian} Lie algebra.
\smallskip

Let $L$ be the free Lie algebra generated by $X,Y$.
Let $\hat L$ be the graded completion of $L$.
Introduce the adjoint operators $x=\ad X$, $y=\ad Y$, 
 i.e. $xA=[X,A]$, $yA=[Y,A]$ for $A\in\hat L$.
Denote by $\widehat{\bar L}$ the graded completion of
 the metabelian quotient $\bar L=L/\bigl[[L,L],[L,L]\bigr]$.
\smallskip

\begin{theorem}
Under $L\to\widehat{\bar L}$ 
 the Hausdorff series $H=\ln(e^Xe^Y)$ 
$$\mbox{maps onto }\quad\bar H=X+Y+\frac{1}{y}
 \left(1-\frac{e^x-1}{x}\cdot\frac{x+y}{e^{x+y}-1}\right)[XY],
 \quad\mbox{ where}$$
 the operator acting on $[XY]$ is considered as
 a commutative series in $x,y$.
\end{theorem}
\smallskip

The series in $x,y$ is a genuine power series with non-negative powers.
The key advantage of the above \emph{metabelian} BCH formula
 is that $H$ is written in a linear basis of
 the free metabelian Lie algebra $\bar L$, see Claim~3.1.
The metabelian BCH formula can be effectively applied
 for solving exponential equations in Lie algebras,
 see Propositions~4.2, 4.4.
\smallskip

In section 2 we recall classical versions of the BCH formula.
The original proof of Theorem~1.1 used large combinatorial formulae 
 involving extended Bernoulli numbers, see
 \cite[Propositions~2.8, 2.12]{Kur}.
Section~3 contains a simpler proof of Theorem~1.1.
In section 4 we give applications to 
 solving exponential equations in metabelian Lie algebras.
\smallskip

\noindent
The author thanks A.~Alekseev, S.~Garoufalidis, H.~Morton for useful remarks.
The author was supported by Marie Curie Fellowship 007477.


\section{The classical Baker-Campbell-Hausdorff formula}


\subsection{The recursive Baker-Campbell-Hausdorff formula}
\noindent
\smallskip

Here we recall the classical version of the BCH formula
 (Theorem~2.4) originally proved by H.~Baker \cite{Bak},
 J.~Campbell \cite{Cam} and F.~Hausdorff \cite{Hau}.

\begin{definition}
The \emph{Bernoulli} numbers $B_n$ are defined by
 the generating function:
 $\sum\limits_{n=0}^{\infty}\dfrac{B_n}{n!}t^n=\dfrac{t}{e^t-1}$,
 e.g. $B_0=1$, $B_1=-\frac{1}{2}$, $B_2=\frac{1}{6}$, $B_3=0$.
\end{definition}

One can verify that $\dfrac{t}{e^t-1}+\dfrac{t}{2}$ is
 an even function, hence $B_n=0$ for all odd $n\geq 3$.
The Bernoulli numbers can be easily computed from
 the recursive relation
 $\sum\limits_{n=1}^m\binom{m+1}{n} B_n=-1$, $m\geq 1$,
 see \cite[Lemma~2.2a]{Kur}.

\begin{definition}
A \emph{derivation} of a Lie algebra $L$ is a linear function
 $D:\hat L\to\hat L$ satisfying the Leibnitz rule
 $D([A,B])=[D(A),B]+[A,D(B)]$.
\end{definition}

The \emph{adjoint} operator
 $\ad A:\hat L\to\hat L$, $\ad A(B)=[A,B]$, is a derivation.
The Leibnitz rule for the derivation $\ad A$ coincides with
 the Jacobi identity: $\ad A([B,C])=[\ad A(B),C]+[B,\ad A(C)]$.
Any derivation of a free Lie algebra can be defined by 
 its values on generators.

\begin{definition}
Let $L$ be the free Lie algebra generated by $X$ and $Y$.
Denote by $D_Y$ the derivation of $L$
 \emph{with respect to} $Y$ such that $D_Y(X)=0$ and 
 $D_Y(Y)=H_1=X+\sum\limits_{n=1}^{\infty}\dfrac{B_{n}}{n!}[Y^nX]\in\hat L$.
\end{definition}

Due to the property $B_n=0$ for odd $n\geq 3$,
 the series $H_1$ can be rewritten as follows:
 $D_Y(Y)=H_1=X+\dfrac{1}{2}[XY]
 +\sum\limits_{n=1}^{\infty}\dfrac{B_{2n}}{(2n)!}[Y^{2n}X]$.
The following theorem is quoted from
 \cite[Corollaries 3.24--3.25, p.~77--79]{Reu}.

\begin{theorem}
 \cite{Bak, Cam, Hau}
The Hausdorff series $H=\ln(e^Xe^Y)$ is equal to
$$H=\sum\limits_{m=0}^{\infty} H_m, \mbox{ where }
  H_0=Y\mbox{ and } H_m=\dfrac{1}{m} D_Y(H_{m-1})
  \mbox{ for } m\geq 1.\eqno{\qed}$$
\end{theorem}
\smallskip


\subsection{The BCH formulae in the forms of Dynkin and Goldberg}
\smallskip

\begin{theorem}
\cite{Dyn}
The Hausdorff series $H=\ln(e^Xe^Y)$ equals
$$H=\sum\limits_{m=1}^{\infty}\dfrac{(-1)^{m-1}}{m}\sum\limits_{
 \tiny \begin{array}{c} p_i+q_i>0 \\  p_i,q_i\geq 0 \end{array} }
 \dfrac{[X^{p_1}Y^{q_1}X^{p_2}Y^{q_2}\ldots X^{p_m}Y^{q_m}]}
 {p_1!q_1!\ldots p_m!q_m!\sum_{i=1}^m(p_i+q_i)}.
 \eqno{\qed}$$
\end{theorem}
\smallskip

The above Dynkin series is not written in
 a linear basis of the free Lie algebra generated by $X,Y$,
 e.g. $[XY]=-[X^0YXY^0]$.
\smallskip

Let us express the Hausdorff series $H=\ln(e^Xe^Y)$
 via associative monomials in $X,Y$,
 where the exponents $s_1,\dots,s_m$ are positive:
$$H=\sum c_x(s_1,\dots,s_m)X^{s_1}Y^{s_2}\cdots(X\vee Y)^{s_m}+$$
$$ +\sum c_y(s_1,\dots,s_m)Y^{s_1}X^{s_2}\cdots(Y\vee X)^{s_m},\mbox{ where}$$
 $(X\vee Y)^{s_m}$ is $X^{s_m}$ for odd $m$
 and $Y^{s_m}$ for even $m$.
Put $m'=[m/2]$.
\smallskip

\begin{theorem}
\cite{Gol}
The generating function for the coefficients $c_x$ is
$$\sum\limits_{s_1,\dots,s_m} c_x(s_1,\dots,s_m)
  Z_1^{s_1}\cdots Z_m^{s_m}=
  \sum_{i=1}^m Z_i e^{m' Z_i}\prod_{j\neq i}
  \frac{e^{Z_j}-1}{e^{Z_i}-e^{Z_j}}.$$
One has $c_y(s_1,\dots,s_m)=(-1)^{n-1}c_x(s_1,\dots,s_m)$,
 where $n=\sum_{i=1}^m s_i$.
\qed
\end{theorem}
\smallskip

For example, $\sum_{s=1}^{\infty}c_x(s)X^s=X$ and
 $\sum_{s=1}^{\infty}c_y(s)Y^s=Y$, hence the Hausdorff series 
 starts as expected: $H=X+Y+\cdots$.
Theorem~2.6 implies that
 $c_x(s_1,\dots,s_m)$ is invariant under any permutation
 of the $s_i$ and $c_x(s_1,\dots,s_m)=0$ when $m$ is odd
 and $n$ is even.


\section{The metabelian Baker-Campbell-Hausdorff formula}

Here we give an elementary proof of Theorem~1.1.
After we interpret the resulting metabelian BCH formula
 as a linear part of a deeper formula via 
 commutators of commutators, see Proposition~3.8.


\subsection{An elementary proof of the metabelian BCH formula}
\noindent
\smallskip
 
Firstly we describe a linear basis of
 the metabelian quotient $\bar L=L/\big[[L,L],[L,L]\bigr]$
 of the free Lie algebra $L$ generated by $X,Y$.
For any word $W$ of length $\geq 2$, the Jacobi identity implies
 $\bigl[X,[Y,[W]]\bigr]-\bigl[Y,[X,[W]]\bigr]=\bigl[[XY],[W]\bigr]=0$.
So we can permute the letters $X,Y$ and express any long commutator
 via the elements $[X^kY^lXY]$ for $k,l\geq 0$.
Claim~3.1 states that these elements are linearly independent 
 in the quotient $\bar L$.

\begin{claim}
\cite[Theorem~5.7]{Reu}, \cite[Section~4.7]{Bah}\\
Let $L$ be the free Lie algebra generated by $X,Y$.
The metabelian quotient $\bar L=L/\bigl[[L,L],[L,L]\bigr]$
 has the linear basis $X,Y,[X^kY^lXY]$, $k,l\geq 0$.
\qed
\end{claim}

Let us express the Hausdorff series $H=\ln(e^Xe^Y)$ 
 via associative monomials in $X,Y$ as before Theorem~2.6.
Look at the terms of the type $c_{rs}X^rY^s$, where
 the letter $X$ always precedes $Y$.
The generating function $c(u,v)=\sum\limits_{r,s\geq 1}c_{rs}u^rv^s$
 of the corresponding coefficients can be extracted from Theorem~2.6,
 but we prefer to give a simple independent proof.

\begin{lemma}
The generating function $c(u,v)$
 of Golberg's coefficients $c_{rs}$ in front of 
 the terms $X^rY^s$ in the Hausdorff series $\ln(e^Xe^Y)$ 
 is equal to
 $c(u,v)=\sum\limits_{r,s\geq 1} c_{rs}u^rv^s=
 ue^u\dfrac{e^v-1}{e^u-e^v}+ve^v\dfrac{e^u-1}{e^v-e^u}$.
\end{lemma}
\noindent
\emph{Proof.}
Put $P=e^X-1$ and $Q=e^Y-1$.
The `$X$ before $Y$' part of 
$$H=\ln(e^Xe^Y)=\ln(1+P+Q+PQ)=\sum\limits_{m=1}^{\infty} \dfrac{(-1)^{m-1}}{m}(P+Q+PQ)^m$$ 
$$\mbox{equals }\quad \sum\limits_{m=1}^{\infty} \dfrac{(-1)^{m-1}}{m}
 \Bigl(\sum\limits_{\tiny \begin{array}{c} r+s=m\\ r,s\geq 1 \end{array} } P^rQ^s
    +\sum\limits_{\tiny \begin{array}{c} r+s=m-1\\ r,s\geq 0 \end{array} } 
 P^{r+1}Q^{s+1}\Bigr).\eqno{(*)}$$
To simplify the above sum we can consider $P,Q$ as commuting variables.
Claim~3.3 is a straightforward computation,
 multipy both sides by $P-Q$.

\begin{claim}
For commuting variables $P,Q$, the following identities hold: \\
$$\sum\limits_{m=2}^{\infty}\dfrac{(-1)^{m-1}}{m}
 \sum\limits_{\tiny  \begin{array}{c} r+s=m\\ r,s\geq 1 \end{array} } P^r Q^s
=\dfrac{Q\ln(1+P)-P\ln(1+Q)}{P-Q},\leqno{(a)}$$
$$\sum\limits_{m=1}^{\infty}\dfrac{(-1)^{m-1}}{m}
 \sum\limits_{\tiny  \begin{array}{c} r+s=m+1\\ r,s\geq 1 \end{array} } P^r Q^s
=\dfrac{PQ}{P-Q}\ln\dfrac{1+P}{1+Q}.\qquad\square\leqno{(b)}$$
\end{claim}

In $(*)$ replace $P,Q$ by $e^u-1,e^v-1$, respectively,
 where $u,v$ are commuting variables of 
 the generating function $(*)$. 
By Claim~3.3 we have
$$c(u,v)=\dfrac{Q\ln(1+P)-P\ln(1+Q)}{P-Q}+\dfrac{PQ}{P-Q}\ln\dfrac{1+P}{1+Q}=$$
$$=\dfrac{u(e^v-1)-v(e^u-1)+(u-v)(e^u-1)(e^v-1)}{e^u-e^v}=$$
$$=\dfrac{ue^ue^v-ve^ue^v-ue^u+ve^v}{e^u-e^v}
  =ue^u\dfrac{e^v-1}{e^u-e^v}+ve^v\dfrac{e^u-1}{e^v-e^u}
 \mbox{ as required.}\eqno{\square}$$
\medskip

\begin{lemma}
Let $L$ be the free Lie algebra generated by $X,Y$.
The Hausdorff series has the form
 $\ln(e^Xe^Y)=X+Y+\sum\limits_{k,l\geq 0}h_{kl}[X^kY^lXY]+H'$, 
 $h_{kl}\in\Q$, $H'\in[[\hat L,\hat L],[\hat L,\hat L]]$.
The only terms contributing to the monomials
 $c_{rs}X^rY^s$ are $h_{kl}[X^kY^lXY]$.
One has $c_{k+1,l+1}=(-1)^lh_{kl}$ for $k,l\geq 0$.
\end{lemma}
\begin{proof}
The simplest monomials coming from a commutator $[A,B]$
 with $A,B\in L$ have the form $X^rY^s$ or $Y^sX^r$ for $r,s\geq 1$.
The polynomial expression of any commutator
 $\bigl[[A,B],[A',B']\bigr]$ with $A',B'\in L$ starts with monomials
 like $X^rY^{s+s'}X^{r'}$ or $Y^sX^{r+r'}Y^{s'}$.
The result follows from: 
$[X^kY^lXY]=(-1)^lX^{k+1}Y^{l+1}+$(more complicated monomials).
\end{proof}
\medskip

\noindent
{\bf Proof of Theorem~1.1.}
Under $\hat L\to \widehat{\bar L}$ the Hausdorff series 
 maps onto $\bar H=X+Y+\sum\limits_{k,l\geq 0}h_{kl}[X^kY^lXY]$
 for some coefficients $h_{kl}\in\Q$.
By Lemma~3.4 the generating function $c(u,v)$
 of Goldberg's coefficients can be expressed via
 $h(x,y)=\sum\limits_{k,l\geq 0}h_{kl}x^ky^l$
 as $c(x,y)=xyh(x,-y)$.
We consider these commutative series as formal Laurent ones,
 although the results always have non-negative powers.
Lemma~3.2 implies
$$h(x,y)=-\dfrac{1}{xy}c(x,-y)=-\dfrac{1}{xy}
 \left( xe^x\dfrac{e^{-y}-1}{e^x-e^{-y}}
 -ye^{-y}\dfrac{e^x-1}{e^{-y}-e^x}\right)=$$
$$=-\dfrac{e^x}{y}\dfrac{1-e^y}{e^{x+y}-1}
   +\dfrac{1}{x}\dfrac{e^x-1}{1-e^{x+y}}=
  \dfrac{1}{e^{x+y}-1}\left( \dfrac{e^{x+y}-e^x}{y}
  -\dfrac{e^x-1}{x}\right)=$$
$$=\dfrac{1}{y}-\dfrac{e^x-1}{e^{x+y}-1}
 \left(\dfrac{1}{y}+\dfrac{1}{x}\right)=
 \dfrac{1}{y}\left(1-\dfrac{e^x-1}{x}\dfrac{x+y}{e^{x+y}-1}\right)
 \mbox{ as required}. \eqno{\qed}$$

\begin{example}
Theorem~1.1 allows us to compute easily first terms
 of the metabelian series $\bar H$ obtained
 from $\ln(e^Xe^Y)$ by $\hat L\to \widehat{\bar L}$.
Suffice to expand the operator $h(x,y)$ acting
 on $[XY]$ in Theorem~1.1, namely \\
 $\bar H=X+Y+h(x,y)[XY]$, where $h(x,y)=$
$$\dfrac{1}{2}+\dfrac{x-y}{12}-\dfrac{xy}{24}
 -\dfrac{x^3+4x^2y-4xy^2-y^3}{720}
 +\dfrac{x^3y+4x^2y^2+xy^3}{1440}
 +(\mbox{deg}\geq 5).$$

\noindent
The above terms agree with the results from
 \cite[Appendix, Proposition A.4]{Kur} and
 coincide with the commutators of $\ln(e^Xe^Y)$ up to degree~4.
\end{example}


\subsection{A deeper BCH formula via commutators of commutators}
\noindent
\smallskip

The original BCH formula of Theorem~2.4 gives a recursive way
 to compute $\ln(e^Xe^Y)$ via commutators in $X,Y$.
The linear part of this formula is $X+Y$.
Here we interpret the metabelian BCH formula of Theorem~1.1
 as a linear part of a deeper formula for $\ln(e^Xe^Y)-X-Y$ 
 via commutators of commutators, see Proposition~3.8.
\smallskip

Let $L$ be the free Lie algebra generated by $X,Y$.
Denote by $\{m,n\}$ the long commutator $[X^mY^{n+1}X]$,
 e.g. $\{0,0\}=[YX]$.
The series $H_1$ of Theorem~2.4 can be rewritten as follows:
 $H_1=X+\sum\limits_{n=1}^{\infty}\dfrac{B_n}{n!}\{0,n-1\}$.
Lemma~3.6 computes the action of $x,y$ on long commutators.

\begin{lemma}
The operators $x,y$ act as follows: $x\{m,n\}=\{m+1,n\}$,\\ 
$y\{m,n\}=\{m,n+1\}+\sum\limits_{k=1}^m\binom{m}{k}[\{k-1,0\},\{m-k,n\}]$. 
\end{lemma}
\begin{proof}
The first formula is trivial, the second one is obtained by induction on $m$.
The base: $y\{0,n\}=[Y^{n+2}X]=\{0,n+1\}$.
The inductive step: 
\smallskip

\noindent
$y\{m+1,n\}=[YX^{m+1}Y^{n+1}X]=\bigl[[YX],[X^mY^{n+1}X]\bigr]+xy[X^mY^{n+1}X]$
\smallskip

\noindent
$=[\{0,0\},\{m,n\}]+x\{m,n+1\}+\sum\limits_{k=1}^m\binom{m}{k}x[\{k-1,0\},\{m-k,n\}]=$
\smallskip

\noindent
$[\{0,0\},\{m,n\}]+\sum\limits_{k=1}^{m}\binom{m}{k}
 \bigl([\{k,0\},\{m-k,n\}]+[\{k-1,0\},\{m-k+1,n\}]\bigr)+$
\smallskip

\noindent
$+\{m+1,n+1\}=\sum\limits_{k=1}^{m+1}
 \Bigl(\binom{m}{k-1}+\binom{m}{k}\Bigr)[\{k-1,0\},\{m-k+1,n\}]+$
\smallskip

\noindent
$+\{m+1,n+1\}=\{m+1,n+1\}+\sum\limits_{k=1}^{m+1}\binom{m+1}{k}[\{k-1,0\},\{m-k+1,n\}]$.
\end{proof}

Lemma~3.6 allows us to compute $x^ky^l\{m,n\}$ for all $k,l>0$
 by using the Leibnitz rule, e.g. 
 $y[\{m,n\},\{r,s\}]=[y\{m,n\},\{r,s\}]+[\{m,n\},y\{r,s\}]$.
We rewrite the derivative $D_Y$ via the operators
 $x=\ad X$, $y=\ad Y$. 

\begin{lemma}
The derivative $D_Y$ mapping $X$ and $Y$ 
 to 0 and $H_1$, respectively, can be expressed via 
 the adjoint operators $x,y$ as follows:
\smallskip

\noindent
 $D_Y\{m,0\}= -\sum\limits_{l=1}^{\infty}\dfrac{B_l}{l!}\{m+1,l-1\},
 \; D_Y\{m,n\}=-x^my^n\sum\limits_{l=1}^{\infty}\dfrac{B_l}{l!}\{1,l-1\}+$
\smallskip

\noindent
$+x^m\sum\limits_{k=0}^{n-1}y^k\Bigl(\{1,n-k-1\}+
 \sum\limits_{l=1}^{\infty}\dfrac{B_l}{l!}\bigl[\{0,l-1\},\{0,n-k-1\}\bigr]\Bigr)$
 for $n\geq 1$.
\end{lemma}
\begin{proof}
We have $D_Y([X^mYX])=-[X^{m+1}H_1]=
 -x^{m+1}\sum\limits_{l=1}^{\infty}\dfrac{B_l}{l!}\{0,l-1\}$.
It remains to apply Lemma~3.6.
The case $n\geq 1$ follows from the Leibnitz rule:
 $D_Y([X^mY^{n+1}X])=x^my^n[H_1,X]+x^m\sum\limits_{k=0}^{n-1}y^k[H_1,[Y^{n-k-1}X]]$.
\end{proof}

Now we rewrite $\ln(e^Xe^Y)$
 via long commutators and commutators of commutators.
Proposition~3.8 follows directly from Theorem~2.4.

\begin{proposition}
In the graded completion of the free Lie algebra generated by $X,Y$ 
we have
$\ln(e^Xe^Y)=X+Y+\sum\limits_{n=1}^{\infty}\dfrac{B_n}{n!}\{0,n-1\}
 +\sum\limits_{m=2}^{\infty}H_{(m)}$, where
 $H_{(2)}=\dfrac{1}{2}\sum\limits_{k=1}^{\infty}\dfrac{B_k}{k!}D_Y\{0,k-1\}$,
 $H_{(m+1)}=\dfrac{D_YH_{(m)}}{m+1}$ for $m\geq 2$.
\qed
\end{proposition}

Theorem~1.1 states that the metabelian image of $\ln(e^Xe^Y)$
 is equal to $X+Y-\sum_{kl}h_{kl}\{k,l\}$, where
 the generating function of the coefficients $h_{kl}$ is
 $\sum_{kl}h_{kl}x^ky^l=\dfrac{1}{y}\left(1-\dfrac{e^x-1}{x}\dfrac{x+y}{e^{x+y}-1}\right)$.
This metabelian part can be interpreted
 as an infinite linear combination of long commutators.
All non-linear terms in the formula of Proposition~3.8
 can be rewritten via commutators of long commutators
 due to Lemmas 3.6 and 3.7 
\smallskip

Proposition~3.8 and Lemmas 3.6, 3.7 give a hope to extend
 Theorem~1.1 to the quotient
 $\ti L=L/\bigl[L',[L',L']\bigr]$, $L'=[L,L]$.
A linear basis of the graded completion of $\ti L$ consists of
 $X,Y$, $\{m,n\}$, $[\{k,l\},\{m,n\}]$,
 where either $k>m\geq 0$, $l,n\geq 0$ or $k=m$, $l>n\geq 0$,
 see \cite[section 4.1]{Reu}.


\section{Applications of the metabelian BCH formula}

Here we show how powerful the metabelian BCH formula is 
 for solving exponential equations in metabelian Lie algebras.
We can rewrite both sides of a given equation in a linear basis
 and compare coefficients.


\subsection{The metabelian Zassenhaus formula}
\noindent
\smallskip

A linear combination $C$ of commutators in a free Lie algebra
 is called a \emph{homogeneous} Lie element of \emph{degree} $n$ 
 if the length of all commutators in $C$ is $n$.
The standard definition is equivalent to 
 the above one \cite[Section~1.3]{Reu}.
According to W.~Magnus \cite[section IV]{Mag},
 Zassenhaus proved, but didn't publish
 the following remarkable result.
\smallskip

\begin{theorem}
\emph{(Zassenhaus)}
Let $L$ be the free Lie algebra generated by $X,Y$.
Then $L$ has a uniquely determined homogeneous Lie element
 $C_n$ of degree $n$ for each $n=2,3,4,\dots$, satisfying
 the \emph{Zassenhaus equation}
 $e^{X+Y}=e^Xe^Ye^{C_2}e^{C_3}e^{C_4}\cdots$
 in the universal enveloping algebra of $L$.
\end{theorem}
\smallskip

W.~Magnus proved the above theorem and gave a recursive way
 to compute the elements $C_n$ in \cite[section IV]{Mag}.
We describe explicitly the metabelian images $\bar C_n$
 of the elements $C_n$ in the quotient $\bar L$.
\smallskip

\begin{proposition}
Let $L$ be the free Lie algebra generated by $X,Y$.
In the universal enveloping algebra of $\bar L=L/\bigl[[L,L],[L,L]\bigr]$,
 the solution to the \emph{metabelian} Zassenhaus equation
 $e^{X+Y}=e^Xe^Ye^{\bar C_2}e^{\bar C_3}e^{\bar C_4}\cdots$ is 
$$\sum_{n=2}^{\infty}\bar C_n=
  \dfrac{1}{x+y}\cdot\dfrac{e^{-y}-1}{y}\cdot
  \left( 1+\dfrac{e^{-x}-1}{x}\cdot\dfrac{y}{e^y-1} \right)[XY],$$
 where $\bar C_n\in\bar L$ are homogeneous Lie elements
 of degree $n=2,3,4,\dots$
\end{proposition}
\noindent
\emph{Proof.}
Use Theorem~1.1 in the completion of $\bar L$ :
$Z=\ln(e^{-X}e^{X+Y})=$
$$=(-X)+(X+Y)+ \dfrac{1}{x+y}\left(
 1+\dfrac{e^{-x}-1}{x}\dfrac{y}{e^y-1} \right) [-X,X+Y]=$$
$$=Y-\dfrac{1}{x+y}\left(
 1+\dfrac{e^{-x}-1}{x}\dfrac{y}{e^y-1} \right) [XY],
 \mbox{ where } x=\ad X,\; y=\ad Y.$$
$Z-Y$ belongs to $[L,L]$ and, for $z=\ad Z$, the operator $\dfrac{z-y}{e^{z-y}-1}$ 
 acts identically modulo commutators of commutators.
We perform computations for commutative series 
 in the algebra of Laurent series, though the result
 will be a genuine formal series with non-negative powers:
$$\ln(e^{-Y}e^{-X}e^{X+Y})=\ln(e^{-Y}e^Z)=
  -Y+Z+\dfrac{1}{y}\left( 1-\dfrac{e^{-y}-1}{-y}\right)[-Y,Z]=$$
$$\dfrac{-1}{x+y}\left(
 1+\dfrac{e^{-x}-1}{x}\dfrac{y}{e^y-1} -\bigl( 1+\dfrac{e^{-y}-1}{y}\bigr) 
 \bigl( 1+\dfrac{e^{-x}-1}{x}\dfrac{y}{e^y-1} \bigr) \right) [XY]$$
$$=\dfrac{1}{x+y}\cdot \dfrac{e^{-y}-1}{y}\cdot
  \left( 1+\dfrac{e^{-x}-1}{x}\dfrac{y}{e^y-1} \right)[XY].$$
It remains to notice that $\bar C_n$ consist
 of commutators only, hence
$$\ln(e^{-Y}e^{-X}e^{X+Y})=
  \ln\left( \prod_{n=2}^{\infty}e^{\bar C_n} \right)=
  \sum_{n=2}^{\infty}\bar C_n \mbox{ in }
  \widehat{\bar L} \mbox{ as required.} \eqno{\qed}$$

\begin{example}
Proposition~4.2 allows us to calculate effectively
$$\sum_{n=2}^{\infty}\bar C_n=
 \left(-\dfrac{1}{2}+\dfrac{x}{6}+\dfrac{y}{3}
 -\dfrac{x^2}{24}-\dfrac{xy}{8}-\frac{y^2}{8} \right)[XY]
 +(\mbox{degree}\geq 5),\mbox{ hence}$$
$$\bar C_2=-\dfrac{[XY]}{2},
  \bar C_3=\dfrac{[X^2Y]}{6}+\dfrac{[YXY]}{3},
  \bar C_4=-\dfrac{[X^3Y]}{24}-\dfrac{[XYXY]}{8}-\frac{[Y^2XY]}{8}.$$

\noindent
The elements computed above coincide with
 the original $C_2,C_3,C_4$ since the metabelian quotient
 $\bar L$ contains all commutators up to degree~4.
\end{example}


\subsection{The commutator equation in the metabelian quotient}
\noindent
\smallskip

The famous Kashiwara-Vergne conjecture 
 \cite[p.~250, Proposition~5.3]{KV} involves 
 the following \emph{commutator} equation:
 $\ln(e^Xe^Y)-X-Y=[X,F]+[Y,G]$ for unknown $F,G$ in
 the graded completion of the free Lie algebra $L$ 
 generated by $X,Y$.
M.~Kashiwara and M.~Vergne proved the existence in
 their conjecture for a soluble Lie algebra
 \cite[Proposition~0]{KV}.
Recently A.~Alekseev and E.~Meinrenken showed
 the existence of a solution for any Lie algebra \cite{AM}.
The Kashiwara-Vergne solution was obtained from
 a differential equation not leading to a closed formula.
\smallskip

We solve the commutator equation completely in the graded completion
 of the metabelian quotient $\bar L=L/\bigl[[L,L],[L,L]\bigr]$.
Denote by $\bar H$ the image of the Hausdorff series
 $\ln(e^Xe^Y)$ under $\hat L\to \widehat{\bar L}$.
The \emph{commutator} equation
 $\bar H-X-Y=[X,F(X,Y)]+[Y,G(X,Y)]$
 has the \emph{symmetry} $\{F(X,Y),G(X,Y)\}\lra\{(G(-Y,-X),F(-Y,-X)\}$.
So we can restrict our attention to the \emph{symmetrized} equation
 $\bar H-X-Y=[X,F(X,Y)]+[Y,F(-Y,-X)]$ for $F,G\in\widehat{\bar L}$.
Put $x=\ad X, y=\ad Y$.

\begin{proposition}
Any solution to the equation
 $\bar H-X-Y=[X,F(X,Y)]+[Y,F(-Y,-X)]$
 in $\widehat{\bar L}$ is $F(X,Y)=aX+\dfrac{Y}{4}+f(x,y)[XY]$, where 
$$f(x,y)=\dfrac{1}{y(x-y)}-\dfrac{1}{4x}-\dfrac{e^x-1}{x}\cdot
 \dfrac{x+y}{e^{x+y}-1}\cdot\dfrac{(x+y)e^y+3x-y}{4xy(x-y)} +yg(x,y),$$
 $f(x,y)$ is a commutative series with non-negative powers,
 $a$ is a constant, $g(x,y)$ is any genuine series satisfying $g(x,y)=-g(-y,-x)$.
\end{proposition}

Proposition~4.4 will follow from Lemmas 4.5 and 4.6.

\begin{lemma}
Any solution to $\bar H-X-Y=[X,F(X,Y)]+[Y,F(-Y,-X)]$
 has the form $F(X,Y)=aX+\dfrac{Y}{4}+f(x,y)[XY]$,
 where $a$ is a constant and
 $f(x,y~)$ is a commutative series with non-negative powers, satisfying
$$xf(x,y)-yf(-y,-x)=-\dfrac{1}{2}+
 \frac{1}{y}\left(1-\frac{e^x-1}{x}\cdot\frac{x+y}{e^{x+y}-1}\right).
 \leqno{(4.5)}$$
\end{lemma}
\begin{proof}
Claim~3.1 implies that any series $F\in\widehat{\bar L}$ can be
 written as $F=aX+bY+f(x,y)[XY]$ for some constants $a,b$
 and a commutative series $f(x,y)$.
By Theorem~1.1 the symmetrized equation is equivalent to
 $$\frac{1}{y}\left(1-\frac{e^x-1}{x}\cdot\frac{x+y}{e^{x+y}-1}\right)
 =b-(-b)+xf(x,y)-yf(-y,-x).$$
Since the left hand side starts with $\dfrac{1}{2}$,
 we get $b=\dfrac{1}{4}$ and (4.5) holds.
\end{proof}

\begin{lemma}
In the ring of formal series with non-negative powers, 
 any solution to $xf(x,y)=yf(-y,-x)$ has the form
 $f(x,y)=yg(x,y)$, where $g$ is any function 
 verifying the symmetry $g(-y,-x)=-g(x,y)$.
\end{lemma}
\begin{proof}
The given condition implies that $f(x,y)=yg(x,y)$ 
 for a genuine series $g(x,y)$.
The substitution gives $g(-y,-x)=-g(x,y)$.
\end{proof}
\smallskip

\noindent
{\bf Proof of Proposition~4.4.}
Put $h(x,y)=\dfrac{1}{y}\left(
 1-\dfrac{e^x-1}{x}\cdot\dfrac{x+y}{e^{x+y}-1}\right)$.
The function $h(x,y)$ appeared in Theorem~1.1 and
 satisfies the important symmetry $h(x,y)=h(-y,-x)$
 since $\ln(e^Xe^Y)=-\ln(e^{-Y}e^{-X})$.
Equation~(4.5) is
 $xf(x,y)-yf(-y,-x)=h(x,y)-\dfrac{1}{2}$, hence
$$\left\{\begin{array}{l}
 \odd h(x,y)=x\even f(x,y)-y\even f(-y,-x),\\
 \even h(x,y)-\dfrac{1}{2}=x\odd f(x,y)-y\odd f(-y,-x),
 \end{array}\right. \mbox{ where}$$
$$\odd h(x,y)=\dfrac{h(x,y)-h(-x,-y)}{2}=
 \dfrac{1}{y}-\dfrac{e^x-1}{x}\cdot\dfrac{e^y+1}{2y}\cdot
 \dfrac{x+y}{e^{x+y}-1},$$
$$\even h(x,y)=\dfrac{h(x,y)+h(-x,-y)}{2}
 =\dfrac{e^x-1}{x}\cdot\dfrac{e^y-1}{2y}\cdot
 \dfrac{x+y}{e^{x+y}-1}.$$

The property $\odd h(x,x)=0$ implies that
 the series $\odd h(x,y)$ is divisible by $x-y$.
The function $\even f(x,y)=\dfrac{\odd h(x,y)}{x-y}$
 satisfies the first equation in the above system.
Actually, we have\\
 $\even f(-y,-x)=\dfrac{\odd h(-y,-x)}{(-y)-(-x)}
 =\dfrac{\odd h(x,y)}{x-y}=\even f(x,y)$.
\smallskip

The function $\odd f(x,y)=\dfrac{1}{2x}\left(\even h(x,y)-\dfrac{1}{2}\right)$
 verifies the second equation.
The expression in the brackets is divisible by $x$
 due to $\even h(0,y)=\dfrac{1}{2}$.
So $x\odd f(x,y)-y\odd f(-y,-x)=\even h(x,y)-\dfrac{1}{2}$.
\smallskip

The final solution is
$f(x,y)=\dfrac{\odd h(x,y)}{x-y}+\dfrac{\even
 h(x,y)}{2x}-\dfrac{1}{4x}=$\\
$$=\dfrac{1}{y(x-y)}-\dfrac{1}{4x}-\dfrac{e^x-1}{x}\cdot
 \dfrac{x+y}{e^{x+y}-1}\cdot\dfrac{(x+y)e^y+3x-y}{4xy(x-y)}.\eqno{\qed}$$
\smallskip

We may pose the problem to describe all solutions to 
 the symmetrized equation $\ln(e^Xe^Y)-X-Y=[X,F(X,Y)]+[Y,F(-Y,-X)]$
 in the graded completion $\hat L$ of free Lie algebra $L$ generated by $X,Y$.
A closed BCH formula in a linear basis of the graded completion $\hat L$ 
 would help to find explicitly all solutions to the Kashiwara-Vergne conjecture.
\medskip



\begin{thebibliography}{References}

\bibitem{AM}
A.~ALEKSEEV, E.~MEINRENKEN,
On the Kashiwara-Vergne Conjecture,
math.QA/0506499.
\emph{Invent. Math.} {\bf 164} (2006), no.~3, 615-634.


\bibitem{Bah}
Yu.~BAHTURIN,
\emph{Identical Relations in Lie Algebras},
VNU Science Press, 1987.

\bibitem{Bak}
H.~BAKER,
On a Law of Combination of Operators (second paper),
\emph{Proc. London Math. Soc.} (3) {\bf 1} (1898), N 29, 14--32.

\bibitem{Cam}
J.~CAMPBELL,
Alternants and Continuous Groups,
\emph{Proc. London Math. Soc.} (3) {\bf 2(3)} (1905), 24--47.

\bibitem{Dr1}
V.~DRINFELD,
Quasi-Hopf Algebras,
\emph{Leningrad Math. J.} {\bf 1} (1990), 1419--1457.

\bibitem{Dyn}
E.~DYNKIN,
On the Representaion of the Series $\log(e^Xe^Y)$
 with Non-commuting $X$ and $Y$ by Commutators,
\emph{Mat. Sbornik} {\bf 25} (1949), 155--162.

\bibitem{Gol}
K.~GOLDBERG,
The Formal Power Series for $\log(e^Xe^Y)$,
\emph{Duke Math. J.} {\bf 23} (1956), 13--21.

\bibitem{Hau}
F.~HAUSDORFF,
Die Symbolische Exponentialformel in der Gruppentheorie,
\emph{Leipziger Berichte} {\bf 58} (1906), 19--48.


\bibitem{KV}
M.~KASHIWARA, M.~VERGNE,
The Campbell-Hausdorff Formula and Invariant Hyperfunctions,
\emph{Invent. Math.} {\bf 47} (1978), 249--272.

\bibitem{Kur}
V.~KURLIN,
Compressed Drinfeld Associators,
\emph{J. Algebra} {\bf 292} (2005), 184--242,
available at www.geocities.com/vak26.

\bibitem{Mag}
W.~MAGNUS,
On the Exponential Solution of Differential Equations
 for a Linear Operator,
\emph{Commun. Pure Appl. Math.} {\bf 7} (1954), 649--673.

\bibitem{Reu}
C.~REUTENAUER,
\emph{Free Lie Algebras},
London Math. Soc. Monographs (N.S.) v.~7 (1993).




\end{thebibliography}
\end{document}